\documentclass[12pt,twoside]{amsart}
\usepackage{latexsym,amsmath,amsopn,amssymb,amsthm,amsfonts}
\usepackage[T2A]{fontenc}
\usepackage[cp1251]{inputenc}
\usepackage[english]{babel}
\usepackage{graphicx}

\newtheorem{theorem}{Theorem}

\newtheorem{proposition}{Proposition}

\newtheorem{remark}{Remark}

\newtheorem{corollary}{Corollary}

\newtheorem{lemma}{Lemma}

\DeclareMathOperator{\re}{Re}

\DeclareMathOperator{\im}{Im}

\DeclareMathOperator{\Ad}{Ad}

\DeclareMathOperator{\ad}{ad}

\DeclareMathOperator{\SO}{SO}

\DeclareMathOperator{\SU}{SU}

\DeclareMathOperator{\sgn}{sgn}

\DeclareMathOperator{\Id}{Id}

\DeclareMathOperator{\Exp}{Exp}

\DeclareMathOperator{\conj}{conj}

\DeclareMathOperator{\Conj}{Conj}

\DeclareMathOperator{\Cut}{Cut}

\DeclareMathOperator{\SL}{SL}

\DeclareMathOperator{\loc}{loc}

\DeclareMathOperator{\glob}{glob}

\DeclareMathOperator{\Tr}{Tr}

\begin{document}
\vspace{0.5cm}
\title[Lie groups $\SU(2)\times\mathbb{R}$ and $\SO(3)\times\mathbb{R}$]{Geodesics and shortest arcs of some sub-Riemannian metrics on the Lie groups $\SU(2)\times\mathbb{R}$ and $\SO(3)\times\mathbb{R}$ with three-dimensional generating distributions}
\author{I.~A.~Zubareva}
\thanks{The work was carried out within the framework of the State Contract to the IM SD RAS, project FWNF--2022--0003.}
\address{Sobolev Institute of Mathematics,\newline
	4 Koptyug Av., Novosibirsk, 630090, Russia}
\email{i\_gribanova@mail.ru}

\begin{abstract}
We find geodesics, shortest arcs,  cut loci,  first conjugate loci, distances between arbitrary elements 
for some left-invariant sub-Riemannian metrics on the Lie groups
$\SU(2)\times\mathbb{R}$ and $\SO(3)\times\mathbb{R}$.

\end{abstract}

\maketitle

\section*{Introduction}

In \cite{Mub} G.M.~Mubarakzyanov obtained a classification (up to isomorphism) of all four-dimensional real Lie algebras.
In \cite{BR}, using a modified version of this classification,  Biggs and Remsing described all
automorphisms of Lie algebras under consideration and classified all connected Lie groups with these Lie algebras.
According to \cite{BR}, the four-dimensional real Lie algebra $\mathfrak{g}_{3,7}\oplus\mathfrak{g}_1$ has a basis
$E_1$,\,$E_2$,\,$E_3$,\,$E_4$ satisfying the following commutation relations:
\begin{equation}
\label{abc}
[E_1,E_2]=E_3,\quad [E_2,E_3]=E_1,\quad [E_3,E_1]=E_2,\quad [E_i,E_4]=0,\,\,i=1,2,3.
\end{equation}
Obviously, $\mathfrak{g}_{3,7}\oplus \mathfrak{g}_1\cong\mathfrak{su}(2)\oplus\mathfrak{\mathbb{R}}$.

According to \cite{BR}, there are five types of connected Lie groups $G$ with Lie algebra $\mathfrak{su}(2)\oplus\mathfrak{\mathbb{R}}$, namely,

1) the universal covering group $\SU(2)\times\mathbb{R}$;

2) $\SU(2)\times\mathbb{R}/\langle\exp(2\pi E_1)\rangle\cong \SO(3)\times\mathbb{R}$;

3)  $\SU(2)\times\mathbb{R}/\langle\exp(E_4)\rangle\cong \SU(2)\times\mathbb{T}$;

4) $\SU(2)\times\mathbb{R}/\langle\exp(2\pi E_1)\exp(E_4)\rangle\cong{\rm U}(2)$;

5) $\SU(2)\times\mathbb{R}/\langle\exp(2\pi E_1)\rangle\langle\exp(E_4)\rangle\cong \SO(3)\times\mathbb{T}$.

The subspaces $\mathfrak{q}_1$ and $\mathfrak{q}_2$ 2 of the Lie algebra are said to be equivalent if $\mathfrak{q}_1$ can be transformed into $\mathfrak{q}_2$ by some automorphism of the Lie algebra. As proved in
\cite{BerZub1}, there exist two nonequivalent classes of equivalent  three-dimensional subspaces  generating the algebra 
$\mathfrak{su}(2)\oplus\mathfrak{\mathbb{R}}$  by the Lie  bracket $[\cdot,\cdot]$.
Moreover, $\mathfrak{q}$ belongs to the first (the second) equivalence class provided that
$\mathfrak{\mathbb{R}}\not\subset\mathfrak{q}$ (respectively, $\mathfrak{\mathbb{R}}\subset\mathfrak{q}$). 

Put $\mathfrak{q}_1={\rm span}\{E_1,E_4-E_3,E_2\}$, $\mathfrak{q}_2={\rm span}\{E_1,E_4,E_2\}$. It is clear that the subspace $\mathfrak{q}_1$ ($\mathfrak{q}_2$) of the Lie algebra $\mathfrak{su}(2)\oplus\mathfrak{\mathbb{R}}$ belongs to the first (respectively, the second) equivalence class.

Using the results of [4], we found geodesics of sub-Riemannian metrics $d_i$, $i=1,2$, on the Lie group $G$ (with the unit ${\rm Id}$)  with the Lie algebra $\mathfrak{su}(2)\oplus\mathfrak{\mathbb{R}}$ which are 
defined by completely nonholonomic left-invariant distributions
$D_i$ with $D_i(\Id)=\mathfrak{q}_i$ and the inner product
$\langle\cdot,\cdot\rangle_i$ on $D_i(\Id)$ with orthonormal basis $E_1,E_4-E_3,E_2$ for  $i=1$ and $E_1,E_4,E_2$ for $i=2$.
 Applying  the methods and results of [5--8], we found shortest arcs, the cut loci, the first conjugate sets, distances between arbitrary elements on the Lie groups
$SU(2)\times\mathbb{R}$ and $SO(3)\times\mathbb{R}$ equipped with these metrics.

\section{Preliminaries}

\begin{theorem}
\label{class1}
Let the basis
\begin{equation}
\label{abc1}
e_1=E_1,\quad e_2=E_4-E_3,\quad e_3=E_2,\quad e_4=E_3
\end{equation}
be given of the Lie algebra $\mathfrak{g}=\mathfrak{su}(2)\oplus\mathfrak{\mathbb{R}}$, $D_1(\Id)={\rm span}(e_1,e_2,e_3)$,
and the inner product $\langle\cdot,\cdot\rangle_1$ on $D_1(\Id)$ with orthonormal basis $e_1,e_2,e_3$.
Then the left-invariant distribution  $D_1$ on the connected Lie group $G$ (with the unit $\Id$) with the Lie algebra $\mathfrak{g}$ with given
$D_1(\Id)$ is totally nonholonomic and the pair $(D_1(\Id),\langle\cdot,\cdot\rangle_1)$ defines the left-invariant sub-Riemannian metric 
$d_1$ on $G$. Moreover, each geodesic $\gamma_1=\gamma_1(t)$, $t\in\mathbb{R}$, in  $(G,d_1)$ parametrized by the arclength such that  
$\gamma_1(0)=\Id$ is a product of two 1--parameter subgroups:
\begin{equation}
\label{geod}
\gamma_1(t)=\gamma_1(\alpha_1,\alpha_2,\alpha_3,\beta;t)=\exp(t(\alpha_1e_1+\alpha_2e_2+\alpha_3e_3+\beta e_4))\exp(-t\beta e_4),
\end{equation}
where $\alpha_1,\alpha_2,\alpha_3$ and $\beta$ are some arbitrary constants such that
\begin{equation}
\label{norm}
\alpha_1^2+\alpha_2^2+\alpha_3^2=1.
\end{equation}
\end{theorem}

\begin{proof}
It follows from (\ref{abc}) and (\ref{abc1}) that	
\begin{equation}
\label{e4}
\left[e_1,e_2\right]=-\left[e_1,e_4\right]=e_3,\quad [e_1,e_3]=e_4,\quad [e_2,e_3]=[e_3,e_4]=e_1,\quad [e_2,e_4]=0.
\end{equation}
This implies the first statement of Theorem 1.	
	
By (\ref{e4}) and \cite[Theorem 3]{BerZub1} every abnormal extremal of the sub-Riemannian space $(G,d_1)$ is nonstrictly abnormal, 
and so  it is a geodesic. According to \cite[Theorem~9]{BerZub2}, every normal geodesic $\gamma_1(t)$, $t\in\mathbb{R}$, with $\gamma_1(0)=\Id$, of $(G,d_1)$ satisfies the differential equation
\begin{equation}
\label{basic}
\dot{\gamma}(t)=dl_{\gamma(t)}(u(t)),\quad u(t)=\psi_1(t)e_1+\psi_2(t)e_2+\psi_3(t)e_3,\quad |u(0)|=1,
\end{equation}
and the absolutely continuous functions $\psi_i(t)$, $i=1,\dots,4$, satisfy the system of differential equations
\begin{equation}
\label{cosystem}
\dot{\psi}_j(t)=\sum_{k=1}^{4}\sum_{i=1}^{3}C_{ij}^k\psi_i(t)\psi_k(t),\quad j=1,\dots,4.
\end{equation}
Here $C_{ij}^k$ are  structure constants in the basis $e_1,e_2,e_3,e_4$ for $\mathfrak{g}$.	

Taking into account (\ref{e4}) the system (\ref{cosystem}) is written as
$$\dot{\psi}_1=-\psi_3(\psi_2+\psi_4),\quad \dot{\psi}_2=0,\quad\dot{\psi}_3=\psi_1(\psi_2+\psi_4),\quad \dot{\psi}_4=0.$$
Assign an arbitrary set of initial data 
$\psi_i(0)=\alpha_i,\,\,i=1,2,3,$ and $\psi_4(0)=\beta-\alpha_2$
of the system.
It is easy to see that $\psi_4(t)\equiv\beta-\alpha_2$,
\begin{equation}
\label{psi123}
\psi_1(t)=\alpha_1\cos\beta t-\alpha_3\sin\beta t,\quad\psi_2(t)\equiv\alpha_2,\quad
\psi_3(t)=\alpha_1\sin\beta t+\alpha_3\cos\beta t,
\end{equation}
and the condition $|u(0)|=1$ is equivalent to (\ref{norm}).

Let us prove that (\ref{geod}) is a solution of (\ref{basic}). It follows from (\ref{e4}) that 
$(\ad(e_4))=e_{31}-e_{13}$, where $(f)$ denotes the matrix of a linear map $f:\mathfrak{g}\rightarrow\mathfrak{g}$ in the base $e_1,e_2,e_3,e_4$; later $(f)$ is identified
with $f$. Here $e_{ij}$ is a $4\times 4$-matrix, having $1$ in the $i$th row and the $j$th column and $0$ in all other.

Differentiating (\ref{geod}) and using (\ref{basic}), (\ref{psi123}), and formula (5) in \cite[Chapter~II, \S~5]{Hel}, we have
$$\dot{\gamma}_1(t)=\exp(t(\alpha_1e_1+\alpha_2e_2+\alpha_3e_3+\beta e_4))(\alpha_1e_1+\alpha_2e_2+\alpha_3e_3+\beta e_4)\exp(-t\beta e_4)+$$ 
$$\gamma_1(t)(-\beta e_4)=\gamma_1(t)\exp(t\beta e_4)(\alpha_1e_1+\alpha_2e_2+\alpha_3e_3+\beta e_4)\exp(-t\beta e_4)+\gamma_1(t)(-\beta e_4)= $$
$$\gamma_1(t)\exp(t\beta e_4)(\alpha_1e_1+\alpha_2e_2+\alpha_3e_3)\exp(-t\beta e_4)+\gamma_1(t)(\beta e_4)+\gamma_1(t)(-\beta e_4)= $$
$$\gamma_1(t)\cdot[\Ad(\exp(t\beta e_4))(\alpha_1e_1+\alpha_2e_2+\alpha_3e_3)]=
\gamma_1(t)\cdot[\exp(t\beta(\ad(e_4)))(\alpha_1e_1+\alpha_2e_2+\alpha_3e_3)]= $$
$$\gamma_1(t)\cdot\left[(\alpha_1\cos\beta t-\alpha_3\sin\beta t)e_1+\alpha_2 e_2+(\alpha_1\sin\beta t+\alpha_3\cos\beta t)e_3\right]=
\gamma_1(t)u(t).$$
\end{proof}

\begin{theorem}
\label{class2}
Let the basis
\begin{equation}
\label{abc2}
e_1=E_1,\quad e_2=E_4,\quad e_3=E_2,\quad e_4=E_3
\end{equation}
be given of the Lie algebra $\mathfrak{g}=\mathfrak{su}(2)\oplus\mathfrak{\mathbb{R}}$, $D_2(\Id)={\rm span}(e_1,e_2,e_3)$,
and the inner product $\langle\cdot,\cdot\rangle_2$  on $D_2(\Id)$ with orthonormal basis $e_1,e_2,e_3$. 
Then the left-invariant distribution  $D_2$ on the connected Lie group $G$ (with the unit $\Id$) with the Lie algebra $\mathfrak{g}$ with given $D_2(\Id)$ is totally nonholonomic and the pair $(D_2(\Id),\langle\cdot,\cdot\rangle_2)$ defines the left-invariant sub-Riemannian metric 
$d_2$ on $G$. Moreover, each geodesic $\gamma_2=\gamma_2(t)$, $t\in\mathbb{R}$, in  $(G,d_2)$ parametrized by the arclength such that  
$\gamma_2(0)=\Id$ is a product of two 1--parameter subgroups:
\begin{equation}
\label{geod2}
\gamma_2(t)=\gamma_2(\alpha_1,\alpha_2,\alpha_3,\beta;t)=\exp(t(\alpha_1e_1+\alpha_2e_2+\alpha_3e_3+\beta e_4))\exp(-t\beta e_4),
\end{equation}
where $\alpha_1,\alpha_2,\alpha_3,\beta$ are some arbitrary constants and (\ref{norm}) holds.
\end{theorem}

\begin{proof}
It follows from (\ref{abc}) and (\ref{abc2}) that	
\begin{equation}
\label{ee4}
\left[e_1,e_2\right]=[e_2,e_3]=[e_2,e_4]=0,\quad [e_1,e_3]=e_4,\quad [e_1,e_4]=-e_3,\quad [e_3,e_4]=e_1.
\end{equation}	
This implies the first statement Theorem 2.		

By (\ref{ee4}) and \cite[Theorem 3]{BerZub1} every abnormal extremal of the sub-Riemannian space $(G,d_2)$ is nonstrictly abnormal and, hence, a geodesic. Moreover, taking into account (\ref{ee4}), the system (\ref{cosystem}) is written as
$$\dot{\psi}_1=-\psi_3\psi_4,\quad \dot{\psi}_2=0,\quad\dot{\psi}_3=\psi_1\psi_4,\quad \dot{\psi}_4=0.$$
Assign an arbitrary set of initial data 
$\psi_i(0)=\alpha_i,\,\,i=1,2,3,$ and $\psi_4(0)=\beta$ of the system.
It is easy to see that $\psi_4(t)\equiv\beta$, and the functions $\psi_i(t)$, $i=1,2,3$, are defined by (\ref{psi123}).
By (\ref{ee4}), $(\ad(e_4))=e_{31}-e_{13}$. To repeat the calculations of the last paragraph of the proof of Theorem \ref{class1}, we obtain that (\ref{geod2}) is the solution of  (\ref{basic}). This completes the proof of Theorem \ref{class2}.
\end{proof}

\begin{proposition}
\label{form}
For all $(\alpha_1,\alpha_2,\alpha_3)\in\mathbb{S}^2$, $\beta\in\mathbb{R}$, and $t\in\mathbb{R}$,
$$\gamma_1(\alpha_1,\alpha_2,\alpha_3,\beta;t)=\gamma_2(\alpha_1,\alpha_2,\alpha_3,\beta-\alpha_2;t)\exp(-t\alpha_2e_4).$$
\end{proposition}

\begin{proof}
By (\ref{abc1}), (\ref{geod}), (\ref{abc2}), and (\ref{geod2}), 
$\gamma_2(\alpha_1,\alpha_2,\alpha_3,\beta-\alpha_2;t)\exp(-t\alpha_2e_4)$
$$=\exp(t(\alpha_1E_1+\alpha_2E_4+\alpha_3E_2+(\beta-\alpha_2)E_3))
\exp(-t(\beta-\alpha_2)E_3)\exp(-t\alpha_2E_3)$$
$$=\exp(t(\alpha_1E_1+\alpha_2(E_4-E_3)+\alpha_3E_2+\beta E_3))\exp(-t\beta E_3)=
\gamma_1(\alpha_1,\alpha_2,\alpha_3,\beta;t).$$
\end{proof}

\begin{proposition}
\label{dop}
Every abnormal extremal in $(G,d_i)$, $i=1,2$, 
parametrized by the arclength is  nonstrictly abnormal and is one of two one-parameter subgroups
$$\gamma_i(0,\alpha_2,0,\beta;t)=\exp(t\alpha_2e_2),\quad\alpha_2=\pm 1,\,\,\,\beta\in\mathbb{R},$$
or their left shifts. 
\end{proposition}

\begin{proof}
	If $\alpha_2=\pm 1$ then  $\alpha_1=\alpha_3=0$ due to (\ref{norm}).  It follows from (\ref{geod}),  (\ref{e4}), (\ref{geod2}), (\ref{ee4}) and the equality 
	$\exp(X+Y)=\exp{X}\exp{Y}$, if $[X,Y]=0$, $X,Y\in\mathfrak{g}$ (see, for example, \cite[Chapter II, $\S$ 5]{Hel}) that for any
	$\beta\in\mathbb{R}$ we have
	$$\gamma_i(0,\pm 1,0,\beta;t)=\exp(t(\pm e_2+\beta e_4))\exp(-t\beta e_4)=\exp(\pm te_2).$$
	This equality and \cite{BerZub1} imply the remaining statements of Proposition \ref{dop}.
\end{proof}

\begin{proposition}
\label{dop2}
Let $\gamma_i(t)=\gamma_i(\alpha_1,\alpha_2,\alpha_3,\beta;t)$, $t\in\mathbb{R}$, be a geodesic in $(G,d_i)$, $i=1,2$, defined by (\ref{geod}) and (\ref{geod2}). 
Then for every $t_0\in\mathbb{R}$ we have
$$\gamma_i(t_0)^{-1}\gamma_i(t)=\gamma_i(\alpha_1\cos\beta t_0-\alpha_3\sin\beta t_0,\alpha_2,\alpha_1\sin\beta t_0+\alpha_3\cos\beta t_0,\beta;t-t_0).$$
\end{proposition}

\begin{proof}
Using (\ref{geod}), (\ref{geod2}), and the equality $(\ad(e_4))=e_{31}-e_{13}$ in the proofs of Theorems \ref{class1}, \ref{class2}, and  (5) in \cite[Chapter II, $\S$ 5]{Hel}, we get
$$\gamma_i(t_0)^{-1}\gamma_i(t)=\exp(t_0\beta e_4)\exp(-t_0(\alpha_1e_1+\alpha_2e_2+\alpha_3e_3+\beta e_4))\times$$
$$\exp(t(\alpha_1e_1+\alpha_2e_2+\alpha_3e_3+\beta e_4))\exp(-t\beta e_4)=$$
$$\exp(t_0\beta e_4)\exp((t-t_0)(\alpha_1e_1+\alpha_2e_2+\alpha_3e_3+\beta e_4))\exp(-t_0\beta e_4)\exp(-(t-t_0)\beta e_4)=$$
$$\exp[\Ad(\exp(t_0\beta e_4))((t-t_0)(\alpha_1e_1+\alpha_2e_2+\alpha_3e_3+\beta e_4))]\cdot\exp(-(t-t_0)\beta e_4)=$$
$$\exp[\exp(\ad(t_0\beta e_4))((t-t_0)(\alpha_1e_1+\alpha_2e_2+\alpha_3e_3+\beta e_4))]\cdot\exp(-(t-t_0)\beta e_4)=$$
$$\exp((t-t_0)((\alpha_1\cos\beta t_0-\alpha_3\sin\beta t_0)e_1+\alpha_2e_2+(\alpha_1\sin\beta t_0+\alpha_3\cos\beta t_0)e_3+\beta e_4)\cdot$$
$$\exp(-(t-t_0)\beta e_4)=\gamma_i(\alpha_1\cos\beta t_0-\alpha_3\sin\beta t_0,\alpha_2,\alpha_1\sin\beta t_0+\alpha_3\cos\beta t_0,\beta;t-t_0).$$
\end{proof}

\section{Sub-Riemannnian geodesics on $\SU(2)\times\mathbb{R}$}

Recall that the Lie group $\SU(2)$ is a compact simply connected group of all unitary unimodular complex $2\times 2$-matrices:
$$\SU(2)=\left\{(A,B):=\left(\begin{array}{cc}
A & B \\
-\overline{B} & \overline{A}
\end{array}\right)\left|\,\,A,\,B\in\mathbb{C},\,\,|A|^2+|B|^2=1\right.\right\}.$$
The Lie algebra $\mathfrak{su}(2)$ of $\SU(2)$  consists of all skew-hermitian complex $2\times 2$-matrices with zero trace: 
$$\mathfrak{su}(2)=\left\{(X,Y):=\left(\begin{array}{cc}
iX & Y \\
-\overline{Y} & -iX
\end{array}\right)\left|\,\,X\in\mathbb{R},\,\,Y\in\mathbb{C}\right.\right\}.$$
The trivial abelian extension
\begin{equation}
\label{defSU(2)xR}
\SU(2)\times\mathbb{R}=\left\{(A,B,v):=\left(\begin{array}{ccc}
A & B & 0 \\
-\overline{B} & \overline{A} & 0 \\
0 & 0 & e^v
\end{array}\right)\left|\,A,\,B\in\mathbb{C},\,\,|A|^2+|B|^2=1;\,\,v\in\mathbb{R}\right.\right\}
\end{equation}
of $\SU(2)$ is a simply connected Lie group, and its Lie algebra $\mathfrak{su}(2)\oplus\mathfrak{\mathbb{R}}$ has some basis $E_1,E_2,E_3,E_4$ satisfying (\ref{abc}):
\begin{equation}
\label{Basis1}
E_1=\frac{1}{2}(e_{12}-e_{21}),\quad E_2=\frac{i}{2}(e_{12}+e_{21}),\quad E_3=\frac{i}{2}(e_{11}-e_{22}),\quad E_4=e_{33}.
\end{equation}
Here  $e_{ij}$, $i,j=1,\dots,3$, denotes the $3\times 3$-matrix, 
having $1$ in the $i$th row and the $j$th column while all
other entries are $0$.

\begin{remark}
The Lie group $\SU(2)\times\mathbb{R}$ is isomorphic to the multiplicative group of nonzero quaternions.	
\end{remark}

To find geodesics $\tilde{\gamma}_i(t)$, $t\in\mathbb{R}$, in $(\SU(2)\times\mathbb{R},d_i)$, $i=1,2$, 
we need 

\begin{lemma}
\label{help1}
Let $z=(X,Y)\in\mathfrak{su}(2)$, $w:=\sqrt{X^2+|Y|^2}\neq 0$, and let $e$ be the identity  $2\times 2$-matrix. Then
$$\exp(z)=\cos{w}\cdot e+\frac{\sin{w}}{w}z.$$
\end{lemma}

\begin{proof}
The characteristic polynomial of $z$ is equal to
$$P(\lambda)=|z-\lambda e|=\left|\begin{array}{cc}
iX-\lambda & Y \\
-\overline{Y} & -iX-\lambda 
\end{array}\right|=\lambda^2+w^2.$$
By the Hamilton–Cayley Theorem, the matrix $z$ is a root of $P(\lambda)$, i.e., $z^2=-w^2e$. Then
$z^{2n+1}=(-1)^nw^{2n} z$, $z^{2n}=(-1)^{n}w^{2n}e$,  $n\in\mathbb{N}.$ Therefore,
$$\exp(z)=e+\sum\limits_{n=1}^{\infty}\frac{z^n}{n!}=e\sum\limits_{n=0}^{\infty}\frac{(-1)^{n}w^{2n}}{(2n)!}+
\frac{z}{w}\sum\limits_{n=0}^{\infty}\frac{(-1)^{n}w^{2n+1}}{(2n+1)!}=\cos{w}\cdot e+\frac{\sin{w}}{w}z.$$
\end{proof}

\begin{theorem}
\label{mainn}
Let
\begin{equation}
\label{mn2}
w_2=\sqrt{1-\alpha_2^2+\beta^2},\quad n_2=\cos\frac{w_2t}{2},\quad 
m_2=\frac{1}{w_2}\sin\frac{w_2t}{2}.
\end{equation}
If $\alpha_2\neq\pm 1$ then the geodesic $\tilde{\gamma}_2(t)=\tilde{\gamma}_2(\alpha_1,\alpha_2,\alpha_3,\beta;t)=(A,B,v)(t)$, $t\in\mathbb{R}$,  of the sub-Riemannian space $(\SU(2)\times\mathbb{R},d_2)$ is defined by the formulas
\begin{equation}
\label{fform1}
A=\left(n_2\cos\frac{\beta t}{2}+\beta m_2\sin\frac{\beta t}{2}\right)+
\left(\beta m_2\cos\frac{\beta t}{2}-n_2\sin\frac{\beta t}{2}\right)i,
\end{equation}
\begin{equation}
\label{fform2}
B=m_2\sqrt{1-\alpha_2^2}\left[\cos\left(\frac{\beta t}{2}+\varphi_0\right)+\sin\left(\frac{\beta t}{2}+\varphi_0\right)i\right],\quad v=\alpha_2t,
\end{equation}
where 
\begin{equation}
\label{varphi}
\cos\varphi_0=\alpha_1/\sqrt{1-\alpha_2^2},\quad \sin\varphi_0=\alpha_3/\sqrt{1-\alpha_2^2}.
\end{equation}
	
If $\alpha_2=\pm 1$ then the geodesic  of the sub-Riemannian space $(\SU(2)\times\mathbb{R},d_2)$ is equal to
$\tilde{\gamma}_2(t)=(1,0,\alpha_2 t)$, $t\in\mathbb{R}$.
\end{theorem}

\begin{proof}
Let $\alpha_2\neq \pm 1$. It follows from (\ref{norm}), (\ref{abc2}),  (\ref{Basis1}), (\ref{mn2}) and Lemma \ref{help1} that
$$\exp(t(\alpha_1e_1+\alpha_2e_2+\alpha_3e_3+\beta e_4))=
\left(\begin{array}{ccc}
n_2+\beta m_2i & (\alpha_1+\alpha_3i)m_2 & 0 \\
(-\alpha_1+\alpha_3i)m_2 & n_2-\beta m_2 i & 0 \\
0 & 0 & e^{\alpha_2t}
\end{array}\right),$$
$$\exp(-t\beta e_4)=
\left(\begin{array}{ccc}
\cos\frac{\beta t}{2}-i\sin\frac{\beta t}{2} & 0 & 0 \\
0 & \cos\frac{\beta t}{2}+i\sin\frac{\beta t}{2} & 0 \\
0 & 0 & 1
\end{array}\right).$$
By Theorem \ref{class2}, it remains to multiply these matrices.

In the case $\alpha_2=\pm 1$ we need to use Proposition \ref{dop}, (\ref{abc2}), (\ref{Basis1}) and Lemma \ref{help1}.
\end{proof}

The following is immediate from  Theorem \ref{mainn} and Propositions \ref{form}, \ref{dop}:

\begin{theorem}
\label{main}
Let
\begin{equation}
\label{mn1}
w_1=\sqrt{1-\alpha_2^2+(\beta-\alpha_2)^2},\quad n_1=\cos\frac{w_1t}{2},\quad 
m_1=\frac{1}{w_1}\sin\frac{w_1t}{2}.
\end{equation}
If $\alpha_2\neq\pm 1$ then the geodesic $\tilde{\gamma}_1(t)=\tilde{\gamma}_1(\alpha_1,\alpha_2,\alpha_3,\beta;t)=(A,B,v)(t)$, $t\in\mathbb{R}$,  of the sub-Riemannian space $(\SU(2)\times\mathbb{R},d_1)$ is defined by the formulas
\begin{equation}
\label{form1}
A=\left(n_1\cos\frac{\beta t}{2}+(\beta-\alpha_2)m_1\sin\frac{\beta t}{2}\right)+
\left((\beta-\alpha_2)m_1\cos\frac{\beta t}{2}-
n_1\sin\frac{\beta t}{2}\right)i,
\end{equation}
\begin{equation}
\label{form2}
B=m_1\sqrt{1-\alpha_2^2}\left[\cos\left(\frac{\beta t}{2}+\varphi_0\right)+\sin\left(\frac{\beta t}{2}+\varphi_0\right)i\right],\quad v=\alpha_2t,
\end{equation}
and (\ref{varphi}) holds.

If $\alpha_2=\pm 1$ then the geodesic $\tilde{\gamma}_1(t)=(A,B,v)(t)$, $t\in\mathbb{R}$, of the sub-Riemannian space $(\SU(2)\times\mathbb{R},d_1)$  is defined by the formulas
\begin{equation}
\label{form3}
A=e^{-i(\alpha_2t/2)},\quad B\equiv 0,\quad v=\alpha_2t.
\end{equation}
\end{theorem}

\begin{remark}
\label{zam0}	
In the case $\alpha_2\neq\pm 1$ we will also denote the geodesic
 $\tilde{\gamma}_i(\alpha_1,\alpha_2,\alpha_3,\beta;t)$ by
$\tilde{\gamma}_i(\varphi_0,\alpha_2,\beta;t)$, if (\ref{varphi}) holds.
\end{remark}

\section{The Cut Locus and First Conjugate Locus in $\SU(2)\times\mathbb{R}$}

In this section, we will recall necessary definitions from the paper \cite{Sachkov} by  Sachkov  for the Lie group $G$ considered in Section 2.

The time moment $\widehat{t}>0$ is called the conjugate time for the normal geodesic
$\gamma_i(t)=\Exp_i(\lambda,t)$, $\lambda=(\alpha_1,\alpha_2,\alpha_3,\beta)\in C:=\mathbb{S}^2\times\mathbb{R}$, 
of the sub-Riemannian space $(G,d_i)$, $i=1,2$, 
if $(\lambda,\widehat{t})$ is the critical point of the exponential map, i.e. 
the differential 
$$(\Exp_i)_{\ast(\lambda,\widehat{t})}:T_{(\lambda,\widehat{t})}(C\times\mathbb{R}_{+})\rightarrow T_{\widehat{g}_i}G,\quad\text {where } \widehat{g}_i=\Exp_i(\lambda,\widehat{t}),$$
is degenerate.  The first conjugate time along the normal geodesic $\gamma_i(t)$ is defined by
$$t^1_{\conj}=\inf\{t>0\mid\,t  \text{ is the conjugate time along } \gamma_i(\cdot)\}.$$
If $t^1_{\conj}>0$, then $\gamma_i(t^1_{\conj})$ is called the first conjugate point.

The first conjugate locus of the sub-Riemannian space $(G,d_i)$ is called the set  ${\rm Conj}^{1}_{i}$ of all first
conjugate points along normal geodesics $\gamma_i(t)$ starting at $\Id$.

The cut locus (for $\Id$) of the sub-Riemannian space $(G,d_i)$ is the set $\Cut_i$ of the endpoints $g\in G$ of all shortest arcs
starting at $\Id$ and noncontinuable beyond $g$.

\begin{proposition}
\label{soprtime}
The time moment $\widehat{t}>0$ is the conjugate time along the geodesic
$\tilde{\gamma}_i(t)=\tilde{\gamma}_i(\alpha_1,\alpha_2,\alpha_3,\beta;t)$, $\alpha_2\neq\pm 1$, $t\in\mathbb{R}$, of the sub-Riemannian space $(\SU(2)\times\mathbb{R},d_i)$, 
$i=1,2$, if and only if  
$$\sin\frac{w_i\widehat{t}}{2}\left(\sin\frac{w_i\widehat{t}}{2}-\frac{w_i\widehat{t}}{2}\cos\frac{w_i\widehat{t}}{2}\right)=0,$$
where $w_i$, $i=1,2$, are defined by formulas (\ref{mn2}) and (\ref{mn1}).
\end{proposition}

\begin{proof}
 By Remark 2,  $\tilde{\gamma}_i(t)=\tilde{\gamma}_i(\varphi_0,\alpha_2,\beta;t)$.

Let $i=2$. For convenience of notation, we put
\begin{equation}
\label{AB}
A_1=\re(A),\quad A_2=\im(A),\quad B_1=\re(B),\quad B_2=\im(B)
\end{equation}
and omit the subscript of the functions $\tilde{\gamma}_2(t),w_2(t),n_2(t)$, and $m_2(t)$.

By (\ref{defSU(2)xR}) and (\ref{Basis1}), 
$$\tilde{\gamma}(t)=2B_1(t)E_1+2B_2(t)E_2+2A_2(t)E_3+e^{\alpha_2t}E_4+2A_1(t)E_5,$$
where $E_5=\frac{1}{2}(E-E_4)$ and $E$ is the identity $3\times 3$-matrix.

It follows from (\ref{fform1}) that $(A_1)'_{\varphi_0}=0$ and $(A_2)'_{\varphi_0}=0$.	
It is easy to see that for $t\neq 0$ the vectors $\tilde{\gamma}^{\prime}_{t}$, $\tilde{\gamma}^{\prime}_{\alpha_2}$, $\tilde{\gamma}^{\prime}_{\beta}$,
and $\tilde{\gamma}^{\prime}_{\varphi_0}$ are linearly dependent  if and only if the vectors $t\tilde{\gamma}^{\prime}_{t}-\alpha_2\tilde{\gamma}^{\prime}_{\alpha_2}$,
$\tilde{\gamma}^{\prime}_{\beta}$, and $\tilde{\gamma}^{\prime}_{\varphi_0}$ are linearly dependent. In this case 
\begin{equation}
\label{f}
t\left[(A_1)'_{\beta}(A_2)'_{t}-(A_1)'_{t}(A_2)'_{\beta}\right]-\alpha_2\left[(A_1)'_{\beta}(A_2)'_{\alpha_2}-(A_1)'_{\alpha_2}(A_2)'_{\beta}\right]=0.
\end{equation}
	
By (\ref{mn2}),
$$n^{\prime}_t=-\frac{mw^2}{2},\quad m^{\prime}_t=\frac{n}{2},\quad n^{\prime}_{\alpha_2}=\frac{\alpha_2mt}{2},\quad 
m^{\prime}_{\alpha_2}=\frac{\alpha_2(2m-nt)}{2w^2},$$
$$n^{\prime}_{\beta}=-\frac{\beta mt}{2},\quad m^{\prime}_{\beta}=\frac{\beta(nt-2m)}{2w^2},\quad
n^{\prime}_{\varphi_0}=m^{\prime}_{\varphi_0}=0.$$
Then, using (\ref{fform1}) and (\ref{fform2}), we have
$$(A_1)'_{t}=-\frac{m}{2}(1-\alpha_2^2)\cos\frac{\beta t}{2},\quad
(A_2)'_{t}=\frac{m}{2}(1-\alpha_2^2)\sin\frac{\beta t}{2};$$
$$(A_1)'_{\beta}=\frac{(2m-nt)(1-\alpha_2^2)}{2w^2}\sin\frac{\beta t}{2},\quad
(A_2)'_{\beta}=\frac{(2m-nt)(1-\alpha_2^2)}{2w^2}\cos\frac{\beta t}{2};$$
$$(A_1)'_{\alpha_2}=\frac{\alpha_2 mt}{2}\cos\frac{\beta t}{2}+\frac{\alpha_2\beta(2m-nt)}{2w^2}\sin\frac{\beta t}{2},$$
$$(A_2)'_{\alpha_2}=-\frac{\alpha_2 mt}{2}\sin\frac{\beta t}{2}+\frac{\alpha_2\beta(2m-nt)}{2w^2}\cos\frac{\beta t}{2}.$$
Inserting the found derivatives into  (\ref{f}), we get
$$\frac{m(2m-nt)(1-\alpha_2^2)t}{4w^2}=0,$$
which implies that $m=0$ or $m=nt/2$. In the first case $\tilde{\gamma}^{\prime}_{\varphi_0}=0$. In the second case, it is easy to check that
$$\tilde{\gamma}^{\prime}_{\beta}=mt\sqrt{1-\alpha_2^2}\left[-\sin\left(\frac{\beta t}{2}+\varphi_0\right)E_2+
\cos\left(\frac{\beta t}{2}+\varphi_0
\right)E_3\right]=t\tilde{\gamma}^{\prime}_{\varphi_0}.$$
Thus, in both cases, $\tilde{\gamma}^{\prime}_{t}$, $\tilde{\gamma}^{\prime}_{\alpha_2}$, $\tilde{\gamma}^{\prime}_{\beta}$, and
$\tilde{\gamma}^{\prime}_{\varphi_0}$ are linearly dependent.

The case $i=1$ is considered similarly.
\end{proof}

\begin{corollary}
\label{sopr}
1. The $n$-conjugate time $t^{n}_{\conj}$ along the geodesic $\tilde{\gamma}_i(t)$, $\alpha_2\neq \pm 1$, of the sub-Riemannian space $(\SU(2)\times\mathbb{R},d_i)$, $i=1,2$, has the form
$$t^{2m-1}_{\conj}=\frac{2\pi m}{w_i},\quad t^{2m}_{\conj}=\frac{2x_m}{w_i},\quad m\in\mathbb{N},$$
where $\{x_1,x_2,\dots\}$ are the ascending positive roots of the equation ${\tan}{x}=x$.

2. The first conjugate locus  $\Conj_i^{1}$ of the sub-Riemannian space  $(\SU(2)\times\mathbb{R},d_i)$, $i=1,2$, has the form
$$\Conj_2^{1}=\{(A,B,v)\in\SU(2)\times\mathbb{R}\left|\,|A|=1,\,\,A\neq 1,\,\,v\in\mathbb{R}\right.\},$$
$$\Conj_1^{1}=\{(A,B,v)\in\SU(2)\times\mathbb{R}\left|\,|A|=1,\,\,v\notin -2{\rm Arg}(A)\right.\},$$
where ${\rm Arg}(A)$ is the set of all values of the argument of a complex number $A$.

\end{corollary}

\begin{proof}
The first statement of Corollary  \ref{sopr} follows from Proposition \ref{soprtime}.
Further, by  Theorems \ref{mainn} and \ref{main},
\begin{equation}
\label{conj}
\tilde{\gamma}_j\left(t^1_{\conj}\right)=\tilde{\gamma}_j\left(2\pi/w_j\right)=\left(-e^{-i(\beta\pi/w_{j})},0,2\alpha_2\pi/w_j\right),\quad j=1,2.
\end{equation}

Consider the functions
$$f_j(\alpha_2,\beta)=\frac{\beta\pi}{w_{j}},\quad v_j(\alpha_2,\beta)=\frac{2\alpha_2\pi}{w_{j}},\quad 
(\alpha_2,\beta)\in (-1,1)\times\mathbb{R},\quad j=1,2.$$ 

It follows from (\ref{mn2}) that  ranges of $f_2(\alpha_2,\beta)$ and $v_2(\alpha_2,\beta)$
are the interval $(-\pi,\pi)$ and $\mathbb{R}$ respectively, and
$$v_2(\alpha_2,\beta)=\frac{2\alpha_2\sqrt{\pi^2-f_2^2(\alpha_2,\beta)}}{\sqrt{1-\alpha_2^2}}.$$ 
From here and (\ref{conj}) we obtain $\Conj_2^{1}$.

By (\ref{mn1}) the range of each of the functions
$f_1(\alpha_2,\beta)$ and $v_1(\alpha_2,\beta)$ is the
real line $\mathbb{R}$, and
$$f_1(\alpha_2,\beta)-\frac{1}{2}v_1(\alpha_2,\beta)=
f_2(\alpha_2,\beta-\alpha_2)\in (-\pi,\pi),\quad (\alpha_2,\beta)\in (-1,1)\times\mathbb{R}.$$
From here and the containment $\pi-f_1(\alpha_2,\beta)\in{\rm Arg}(A)$ (see (\ref{conj})) we obtain $\Conj_1^{1}$.

\end{proof}

\begin{proposition}
	\label{ssq3}
	Let $\tilde{\gamma}_1(\varphi_0,\alpha_2,\beta;t)$, $0\leq t\leq T$, be a noncontinuable shortest arc in
	$(SU(2)\times\mathbb{R},d_1)$.
	Then $\tilde{\gamma}_2(\varphi_0,\alpha_2,\beta-\alpha_2;t)$, $0\leq t\leq T$, is a noncontinuable shortest arc in 
	$(SU(2)\times\mathbb{R},d_2)$. The converse is also true.
\end{proposition}

\begin{proof}
	Put $g_1=(A,B,v)=\tilde{\gamma}_1(\varphi_0,\alpha_2,\beta;T)$. By Theorem \ref{main}, $v=\alpha_2T$. Then $g_2:=\tilde{\gamma}_2(\varphi_0,\alpha_2,\beta-\alpha_2;T)=(Ae^{iv/2},Be^{-iv/2},v)$ due to Proposition \ref{form}. 
	Suppose that $\tilde{\gamma}_2(\varphi_0^{\ast},\alpha_2^{\ast},\beta^{\ast};t)$, $0\leq t\leq T_1<T$, is a shortest arc in 	
	$(SU(2)\times\mathbb{R},d_2)$ joining $\Id$ and $g_2$; $v=\alpha_2^{\ast}T_1$ by Theorem \ref{main}. By
	virtue of Proposition \ref{form}, the geodesic segment $\tilde{\gamma}_1(\varphi_0^{\ast},\alpha_2^{\ast},\beta^{\ast}+\alpha_2^{\ast};t)$, 
	$0\leq t\leq T_1$, in $(SU(2)\times\mathbb{R},d_1)$ joins $\Id$ and $g_1$, which is impossible. 
	Hence $\tilde{\gamma}_2(\varphi_0,\alpha_2,\beta-\alpha_2;t)$, $0\leq t\leq T$, is a shortest arc in $(SU(2)\times\mathbb{R},d_2)$. 
	It is easy to check that if  the geodesic segment $\tilde{\gamma}_2(\varphi_0,\alpha_2,\beta-\alpha_2;t)$, $0\leq t\leq T_2$, $T_2>T$, is a 
	shortest arc in $(SU(2)\times\mathbb{R},d_2)$, then $\tilde{\gamma}_1(\varphi_0,\alpha_2,\beta;t)$, $0\leq t\leq T_2$, is a 
	shortest arc in $(SU(2)\times\mathbb{R},d_1)$, which is impossible. Thus $\tilde{\gamma}_2(\varphi_0,\alpha_2,\beta-\alpha_2;t)$, $0\leq t\leq T$,
	is a noncontinuable shortest arc in $(SU(2)\times\mathbb{R},d_2)$. 
	
	The converse statement is proved in a similar way.
\end{proof}	

We need the following proposition from \cite{BZ}:

\begin{proposition}
\label{bz}
If in the Lie group  with left-invariant sub-Riemannian metric 
two points are joined by two different  normal geodesic parametrized by arclength of
equal length, then neither of these
geodesics is a shortest arc nor a part of a longer shortest arc.
\end{proposition}

\begin{proposition}
\label{cut}

1. If $\alpha_2=\pm 1$ then 
the geodesic  
$\tilde{\gamma}_i(t)=\exp(t\alpha_2e_2)$, $t\in\mathbb{R}$, is a metric line of
the sub-Riemannian space $(\SU(2)\times\mathbb{R},d_i)$, $i=1,2$; i.e., the segment 
$\tilde{\gamma}_i(t)=\exp(t\alpha_2e_2)$, $t_0\leq t\leq t_1$, is a shortest arc for all $t_0,t_1\in\mathbb{R}$.

2. Let $\alpha_2\neq\pm 1$ and $\tilde{\gamma}_i=\tilde{\gamma}_i(\varphi_0,\alpha_2,\beta;t)$, $0\leq t\leq T_{i}$, be a noncontinuable shortest arc of
the sub-Riemannian space $(\SU(2)\times\mathbb{R},d_i)$,  $i=1,2$. Then
$$T_{1}=2\pi/\sqrt{1-\alpha_2^2+(\beta-\alpha_2)^2},\quad T_2=2\pi/\sqrt{1-\alpha_2^2+\beta^2}.$$

3. The cut locus $\Cut_i$ in $(\SU(2)\times\mathbb{R},d_i)$, $i=1,2$, coincides with the first conjugate locus $\Conj^{1}_{i}$.
\end{proposition}

\begin{proof}
It follows from Theorems \ref{mainn} and \ref{main} that the length of the arc of the abnormal extremal connecting $\Id$ and $(A,B,v)$, $v\neq 0$, is equal to $|v|$, and the length of the arc of the normal geodesic 
$\tilde{\gamma}_i(t)=\tilde{\gamma}_i(\varphi_0,\alpha_2,\beta;t)$, $\alpha_2\neq\pm 1$, connecting these elements, is equal to $|v/\alpha_2|>|v|$. This implies item~1 of Proposition \ref{cut}.
	
Let $\alpha_2\neq\pm 1$. It follows from (\ref{mn2}) and (\ref{mn1}) that $m_i\left(2\pi/w_i\right)=0$, $i=1,2$; therefore by (\ref{conj}) $\tilde{\gamma}_i\left(2\pi/w_i\right)$ does not depend on $\varphi_0$, and  depends only on $\alpha_2$ and $\beta$. Then  $T_{i}\leq 2\pi/w_i$ by Proposition \ref{bz}.

Let us show that $\tilde{\gamma}_i(t)=\tilde{\gamma}_i(\varphi_0,\alpha_2,\beta;t)$, $0\leq t\leq2\pi/w_i$, $i=1,2$, is a shortest arc. The case $\alpha_2=0$ is considered in \cite{Bosc}, so we assume that $\alpha_2\neq 0$. Put $\tilde{\gamma}_i\left(2\pi/w_i\right)=(A,B,v)$. Then $|A|=1$, $B=0$, and $v\neq 0$. 

Let $i=1$. Inserting $w_1$ (see (\ref{mn1})) into $v=2\alpha_2\pi/w_1$, we get
$$\beta=
\alpha_2+\left(\sgn(\beta-\alpha_2)\sqrt{(4\pi^2+v^2)\alpha_2^2-v^2}\right)/|v|.$$
Then it follows from (\ref{conj}) that $\alpha_2$ satisfies the system of equations
$$\left\{\begin{array}{c}
\cos\left(\frac{v}{2}+\sgn(\beta-\alpha_2)\frac{\sqrt{(4\pi^2+v^2)\alpha_2^2-v^2}}{2|\alpha_2|}\right)=-\re(A), \\
\sin\left(\frac{v}{2}+\sgn(\beta-\alpha_2)\frac{\sqrt{(4\pi^2+v^2)\alpha_2^2-v^2}}{2|\alpha_2|}\right)=\im(A),
\end{array}\right.$$
which is equivalent to the equation
\begin{equation}
\label{equat}
\cos\frac{\sqrt{(4\pi^2+v^2)\alpha_2^2-v^2}}{2|\alpha_2|}=-\cos\left(\arg(A)+v/2\right).
\end{equation}

It is easy to check that the function $f(t)=\frac{1}{2t}\sqrt{(4\pi^2+v^2)t^2-v^2}$ on the interval $\left[|v|/\sqrt{4\pi^2+v^2},1\right]$, is increasing and its range is the segment
$\left[0,\pi\right]$. This implies that (\ref{equat}) has the unique solution $\alpha_2\in (-1,1)$ such that
$\sgn(\alpha_2)=\sgn(v)$. This means that the segment of each geodesic $\tilde{\gamma}_1(\varphi_0,\alpha_2,\beta;t)$  of the sub-Riemannian space
$(\SU(2)\times\mathbb{R},d_1)$ joining $\Id$ and $(A,B,v)$ has the same length $v/\alpha_2$. Hence, it follows from the inequality $T_1\leq 2\pi/w_1$ and Proposition \ref{bz} that each such segment is a shortest arc.

The case $i=2$ is considered similarly.

Item~3 of Proposition \ref{cut} follows from items~1, 2 of Proposition \ref{cut} and Corollary \ref{sopr}.

\end{proof}

\section{Sub-Riemannian Distance on the Lie Group $\SU(2)\times\mathbb{R}$}

In consequence of left invariance of the metrics $d_1$ and $d_2$ it is sufficient to compute the distances $d_1(\Id,g)$ and $d_2(\Id,g)$ between the unit $\Id$ of the Lie group $\SU(2)\times\mathbb{R}$ and an arbitrary element $g=(A,B,v)\in\SU(2)\times\mathbb{R}$.

\begin{proposition}
\label{k}
For every  $(A,B,v)\in\SU(2)\times\mathbb{R}$
$$d_1(\Id,(A,B,v))=d_2(\Id,(Ae^{iv/2},Be^{-iv/2},v)).$$
\end{proposition}

\begin{proof}
Let $(A,B,v)\in\SU(2)\times\mathbb{R}$ and let $\tilde{\gamma}_1(\alpha_1,\alpha_2,\alpha_3,\beta;t)$, $0\leq t\leq t_1$, be a shortest arc in $(\SU(2)\times\mathbb{R},d_1)$ joining $\Id$ and $(A,B,v)$;  and, moreover, $v=\alpha_2 t_1$ due to (\ref{form2}). By (\ref{abc2}) and Proposition \ref{form}, the
geodesic segment $\tilde{\gamma}_2(\alpha_1,\alpha_2,\alpha_3,\beta-\alpha_2;t)$, $0\leq t\leq t_1$, in $(\SU(2)\times\mathbb{R},d_2)$ joins $\Id$ and $(A,B,v)\exp(vE_3)$. By virtue of (\ref{Basis1}), $\exp(vE_3)=(e^{iv/2},0,0)$; therefore $(A,B,v)\exp(vE_3)=(Ae^{iv/2},Be^{-iv/2},v)$. 
We assume that this geodesic segment is not a shortest curve and denote by
$\tilde{\gamma}_2(\alpha_1^{\ast},\alpha_2^{\ast},\alpha_3^{\ast},\beta^{\ast};t)$,
$0\leq t\leq t_1^{\ast}<t_1$, the shortest arc joining $\Id$ and $(A,B,v)\exp(vE_3)$. 
Then it follows from Proposition \ref{form} that the geodesic segment
$\tilde{\gamma}_1(\alpha_1^{\ast},\alpha_2^{\ast},\alpha_3^{\ast},\beta^{\ast}+\alpha_2^{\ast};t)$, $0\leq t\leq t_1^{\ast}$,
in $(\SU(2)\times\mathbb{R},d_1)$ joins $\Id$ and $(A,B,v)$; a contradiction.

\end{proof}

\begin{theorem}
\label{os}
Let $\mathbb{H}\times\mathbb{R}$ be a connected $(n+1)$--dimensional Lie group $($with the unit $\Id=(e,1)$ and the Lie algebra $\mathfrak{g}\oplus\mathfrak{g}_1)$ with the left-invariant sub-Riemannian metric $d$ generated by a completely nonholonomic distribution $D$ with $$D(\Id)={\rm span}(e_1,\dots,e_{n-1},e_{n+1})\subset\mathfrak{g}\oplus\mathfrak{g}_1,\quad e_1,\dots,e_{n-1}\in\mathfrak{g},\,\,e_{n+1}=1\in\mathfrak{g}_1,$$
and the inner product $\langle\cdot,\cdot\rangle$ with the orthonormal basis $e_1,\dots,e_{n-1},e_{n+1}$ on $D(\Id)$.
Then
$d^2(\Id,(g,e^v))=v^2+d^2(\Id,(g,1))$ for all
$(g,e^v)\in\mathbb{H}\times\mathbb{R}$.
\end{theorem}

\begin{proof}
Let $(g_0,e^{v_0})\in\mathbb{H}\times\mathbb{R}$ and $(g_0,e^{v_0})\neq \Id$, while $\gamma(t)=(g,e^v)(t)$, $0\leq t\leq d_0:=d(\Id,(g_0,e^{v_0}))$, is a shortest arc in $(\mathbb{H}\times\mathbb{R},d)$ joining $\Id$ and $(g_0,e^{v_0})$. It follows from \cite[Theorem~9]{BerZub2} that
$\psi_{n+1}(t)\equiv \varphi_{n+1},$ $v(t)=\varphi_{n+1}t$ for some $\varphi_{n+1}\in [-1,1]$, and $g(t)$ is a solution of the system of differential equations
\begin{equation}
\label{gt}
\dot{g}(t)=dl_{g(t)}(u(t)),\quad u(t)=\sum\limits_{i=1}^{n-1}\psi_i(t)e_i,
\end{equation}
and, moreover, the absolutely continuous functions  $\psi_i(t)$, $i=1,\dots,n$, satisfy the Cauchy problem
\begin{equation}
\label{psi}
\dot{\psi}_j(t)=\sum_{k=1}^{n}\sum_{i=1}^{n-1}C_{ij}^k\psi_i(t)\psi_k(t),\quad \psi_j(0)=\varphi_j,\quad j=1,\dots,n;
\end{equation}
here $C_{ij}^k$ are structure constants in the basis  $e_1,\dots,e_n$ of the Lie algebra $\mathfrak{g}$ and
\begin{equation}
\label{psisum}
\sum\limits_{i=1}^{n-1}\psi^2_i(t)=1-\varphi_{n+1}^2.	
\end{equation}	

Since $d(v_0)=d_0$; therefore,   $\varphi_{n+1}=v_0/d_0$. 

If $v_0=0$ then $\varphi_{n+1}=0,$ $v(t)\equiv 0,$ and  $\gamma(t)=(g(t),1)$. Then $g=g(t)$ is a shortest arc in $(\mathbb{H},d_{\mathbb{H}}),$ where $d_{\mathbb{H}}$ is a left-invariant sub-Riemannian metric on the Lie group $\mathbb{H}$ defined by the subspace $D(\Id)\cap\mathfrak{g}\subset\mathfrak{g}$ with the orthonormal basis $e_1,\dots,e_{n-1},$ and $d_\mathbb{H}$ coincides with the metric induced on $\mathbb{H}$ by $d.$ 

If $\varphi_{n+1}=\pm 1$ then $d_0=|v_0|$ and $\psi_i(t)=0$, $i=1,\dots,n-1$, in view of (\ref{psisum}). It follows from (\ref{gt}) that $g(t)$ is a constant and $g(t)\equiv g_0=e$.

In both cases Theorem \ref{os} is true. 

If $0<|\varphi_{n+1}|<1$ then $v_0\neq 0$ and $g_0\neq e$, while, due to  (\ref{gt}) and (\ref{psisum}), the curve $g=g(t)$, $0\leq t\leq d_0$, in $(\mathbb{H},d_\mathbb{H})$ (parameterized in the proportion to the arclength with the factor $l=\sqrt{1-\varphi^2_{n+1}}$) has length $ld_0=\sqrt{d_0^2-v_0^2},$ and besides $g(0)=e$, $g(d_0)=g_0$. 
As a consequence of the above, Theorem \ref{os} is equivalent to the fact that $g=g(t)$, $0\leq t\leq d_0$, is a shortest arc in 
$(\mathbb{H},d_\mathbb{H})$, and this fact is equivalent to the initial assumption that the shortest arc
$\gamma(t)$, $0\leq t\leq d_0$, in $(\mathbb{H}\times\mathbb{R},d)$ has length $d_0\sqrt{l^2+\varphi^2_{n+1}}=d_0$.  
\end{proof}

\begin{corollary}
\label{osn}
$d_2^2(\Id,(A,B,v))=v^2+d_2^2(\Id,(A,B,0))$ for all $(A,B,v)\in\SU(2)\times\mathbb{R}$.
\end{corollary}

The following is immediate from Corollary \ref{osn} and \cite[Theorem~2]{BZ}.

\begin{theorem}
\label{main2}
Let $g=(A,B,v)\in SU(2)\times\mathbb{R}$ and $\arg(A)\in (-\pi,\pi]$. Then

1. If $A=0$ then $d_2(\Id,g)=\sqrt{v^2+\pi^2}$.

2. If $|A|=1$ then $d_2(\Id,g)=\sqrt{v^2+4|\arg(A)|(2\pi-|\arg(A)|)}$.
 
3. If $0<|A|<1$ and $\re(A)=|A|\sin{\left(\frac{\pi}{2}|A|\right)}$ then 
$d_2(\Id,g)=\sqrt{v^2+\pi^2(1-|A|^2)}.$

4. If $0<|A|<1$ and $\re(A)>|A|\sin{\left(\frac{\pi}{2}|A|\right)}$ then
$$d_2(\Id,g)=\sqrt{v^2+\frac{4}{1+\xi^2}\left(\sin^{-1}\sqrt{(1-|A|^2)(1+\xi^2)}\right)^2},$$
where $\xi$ is unique solution of the system of equations
$$\left\{\begin{array}{l}
\cos{\left(-\frac{\xi}{\sqrt{1+\xi^2}}\sin^{-1}\sqrt{(1-|A|^2)(1+\xi^2)}+\sin^{-1}\frac{\xi\sqrt{1-|A|^2}}{|A|}\right)}=\frac{\re(A)}{|A|}, \\
\sin{\left(-\frac{\xi}{\sqrt{1+\xi^2}}\sin^{-1}\sqrt{(1-|A|^2)(1+\xi^2)}+\sin^{-1}\frac{\beta\sqrt{1-|A|^2}}{|A|}\right)}=\frac{\im(A)}{|A|}.
\end{array}\right.$$

5. If $0<|A|<1$ and $\re(A)<|A|\sin{\left(\frac{\pi}{2}|A|\right)}$ then
$$d_2(\Id,g)=\sqrt{v^2+\frac{4}{1+\xi^2}\left(\pi-\sin^{-1}\sqrt{(1-|A|^2)(1+\xi^2)}
\right)^2},$$
where $\xi$ is unique solution of the system of equations
$$\left\{\begin{array}{l}
\cos{\left(\frac{\xi}{\sqrt{1+\xi^2}}\left(\pi-\sin^{-1}\sqrt{(1-|A|^2)(1+\xi^2)}\right)+\sin^{-1}
\frac{\xi\sqrt{1-|A|^2}}{|A|}\right)}=-\frac{\re(A)}{|A|}, \\
\sin{\left(\frac{\xi}{\sqrt{1+\xi^2}}\left(\pi-\sin^{-1}\sqrt{(1-|A|^2)(1+\xi^2)}\right)+\sin^{-1}
\frac{\xi\sqrt{1-|A|^2}}{|A|}\right)}=\frac{\im(A)}{|A|}. 
\end{array}\right.$$
\end{theorem}

Using Proposition \ref{k} and Theorem \ref{main2}, we can obtain exact formulas for $d_1(\Id,(A,B,v))$.
We will not present them because of their bulkiness.

\section{Sub-Riemannian geodesics on $SO(3)\times\mathbb{R}$}

Recall that the Lie group $\SO(3)$ is a connected compact  group of all orthogonal  real $3\times 3$-matrices with 
determinant $1$.
Its Lie algebra $\mathfrak{so}(3)\cong\mathfrak{su}(2)$ consists of all skew-symmetric real $3\times 3$--matrices.

We will deal with the trivial abelian extension of $\SO(3)$:
$$\SO(3)\times\mathbb{R}=\scriptsize{\left\{(C,v):=\left(\begin{array}{ccc}
C & & 0 \\
& & 0 \\
0 & 0 & e^v 
\end{array}\right)\,\,\left|\,\,C\in\SO(3),\,\,v\in\mathbb{R}\right.\right\}.}$$
The Lie algebra  $\mathfrak{so}(3)\oplus\mathbb{R}$ of $\SO(3)\times\mathbb{R}$ has some basis $E_1,E_2,E_3,E_4$ satisfying (\ref{abc}):
\begin{equation}
\label{basis2}
E_1=e_{32}-e_{23},\quad E_2=e_{13}-e_{31},\quad E_3=e_{21}-e_{12},\quad E_4=e_{44}.
\end{equation}
Here $e_{ij}$, $i,j=1,\dots,4$, denotes the $4\times 4$-matrix, 
having $1$ in the $i$th row and the $j$th column and $0$ in all other.

For convenience  we denote by $\rho_i$, $i=1,2$, the left-invariant sub-Riemannian metric on $\SO(3)\times\mathbb{R}$ (with the unit $\Id$)
defined by the completely nonholonomic left-invariant distribution $\Delta_i$ with $\Delta_i(\Id)={\rm span}(e_1,e_2,e_3)$, where the vectors
$e_1,e_2,e_3$ are given by formulas (\ref{abc1}) for $i=1$ and (\ref{abc2}) for $i=2$; we define the inner product $\langle\cdot,\cdot\rangle_i$ on $\Delta_i(\Id)$ with the orthonormal basis $e_1,e_2,e_3$.

To find geodesics in $(\SO(3)\times\mathbb{R},\rho_i)$, $i=1,2$, 
we need 

\begin{lemma}
\label{help}
Let $C=(c_{ij})\in\mathfrak{so}(3)$. Then
$$\exp\left(C\right)=E+\frac{\sin{w}}{w}C+\frac{1-\cos{w}}{w^2}C^2,\quad w=\sqrt{c_{12}^2+c_{13}^2+c_{23}^2},$$
where $E$ is the identity $3\times 3$-matrix.
\end{lemma}

\begin{proof}
The characteristic polynomial of $C$ is equal to
$$P(\lambda)=|C-\lambda E|=\left|\begin{array}{ccc}
-\lambda & c_{12} & c_{13} \\
-c_{12} & -\lambda & c_{23} \\
-c_{13} & -c_{23} & -\lambda
\end{array}\right|=-\lambda^3-w^2\lambda.$$
By the Hamilton–Cayley Theorem,  $C$ is a root of $P(\lambda)$, i.e. $C^3=-w^2C$.  Therefore,
$C^{2n+1}=(-1)^nw^{2n} C$, $C^{2n}=(-1)^{n+1}w^{2n-2}C^2$ for $n\geq 1.$ Then
$$\exp\left(C\right)=E+\sum\limits_{n=1}^{\infty}\frac{C^n}{n!}=E+\frac{C}{w}\sum\limits_{n=0}^{\infty}\frac{(-1)^{n}w^{2n+1}}{(2n+1)!}-
\frac{C^2}{w^2}\sum\limits_{n=1}^{\infty}\frac{(-1)^{n}w^{2n}}{(2n)!},$$
and Lemma \ref{help} is proved.
\end{proof}

\begin{theorem}
\label{mainq}
 Let
 \begin{equation}
 \label{nm2}
 w_2=\sqrt{1-\alpha_2^2+\beta^2},\quad\mu_2=\frac{\sin{w_2t}}{w_2},\quad
 \nu_2=\frac{1-\cos{w_2t}}{w_2^2}.
 \end{equation}
If $\alpha_2\neq\pm 1$ then the geodesic $\gamma_2(t)=\gamma_2(\alpha_1,\alpha_2,\alpha_3,\beta;t)$, $t\in\mathbb{R}$,  of the sub-Riemannian space $(\SO(3)\times\mathbb{R},\rho_2)$  is equal to $(C,v)(t)$, $t\in\mathbb{R}$, where $v(t)=\alpha_2t$ and the columns $C_j(t)$, $j=1,2,3$, of the matrix $C(t)\in\SO(3)$ 
are given by the formulas
$$C_1(t)=\left(\begin{array}{c}
\left(1-\nu_2(\alpha_3^2+\beta^2)\right)\cos{\beta t}-\left(\alpha_1\alpha_3\nu_2-\beta\mu_2\right)\sin{\beta t} \\
\left(\alpha_1\alpha_3\nu_2+\beta\mu_2\right)\cos{\beta t}-\left(1-\nu_2(\alpha_1^2+\beta^2)\right)\sin{\beta t} \\
\left(\alpha_1\beta\nu_2-\alpha_3\mu_2\right)\cos{\beta t}-\left(\alpha_3\beta\nu_2+\alpha_1\mu_2\right)\sin{\beta t} 
\end{array}\right),$$
$$C_2(t)=\left(\begin{array}{c}
\left(1-\nu_2(\alpha_3^2+\beta^2)\right)\sin{\beta t}+\left(\alpha_1\alpha_3\nu_2-\beta\mu_2\right)\cos{\beta t} \\
\left(\alpha_1\alpha_3\nu_2+\beta\mu_2\right)\sin{\beta t}+\left(1-\nu_2(\alpha_1^2+\beta^2)\right)\cos{\beta t} \\
\left(\alpha_1\beta\nu_2-\alpha_3\mu_2\right)\sin{\beta t}+\left(\alpha_3\beta\nu_2+\alpha_1\mu_2\right)\cos{\beta t} 
\end{array}\right),$$
$$C_3(t)=\left(\begin{array}{c}
\alpha_1\beta\nu_2+\alpha_3\mu_2 \\
\alpha_3\beta\nu_2-\alpha_1\mu_2 \\
1-\nu_2(\alpha_1^2+\alpha_3^2) 
\end{array}\right).$$

If $\alpha_2=\pm 1$ then the geodesic of the sub-Riemannian space $(\SO(3)\times\mathbb{R},\rho_2)$ is equal to
$\gamma_2(t)=(E,\alpha_2 t)$, $t\in\mathbb{R}$, where $E$ is the identity $3\times 3$-matrix.	
\end{theorem}

\begin{proof}
Let  $\alpha_2\neq\pm 1$. If follows from (\ref{abc1}), (\ref{basis2}), (\ref{nm2}) and Lemma \ref{help} that
$$\exp(t(\alpha_1e_1+\alpha_2e_2+\alpha_3e_3+\beta e_4))$$
$$=\left(\begin{array}{cccc}
1-\nu_2(\alpha_3^2+\beta^2)  & \alpha_1\alpha_3\nu_2-\beta\mu_2 & \alpha_1\beta\nu_2+\alpha_3\mu_2 & 0 \\
\alpha_1\alpha_3\nu_2+\beta\mu_2 & 1-\nu_2(\alpha_1^2+\beta^2) & \alpha_3\beta\nu_2-\alpha_1\mu_2 & 0 \\
\alpha_1\beta\nu_2-\alpha_3\mu_2 & \alpha_3\beta\nu_2+\alpha_1\mu_2 & 1-\nu_2(\alpha_1^2+\alpha_3^2) & 0 \\
0 & 0 & 0 & e^{\alpha_2 t}
\end{array}\right),$$
$$\exp(-t\beta e_4)=
\left(\begin{array}{cccc}
\cos{\beta t} & \sin{\beta t} & 0 & 0 \\
-\sin{\beta t} & \cos{\beta t} & 0 & 0 \\
0 & 0 & 1 & 0 \\
0 & 0 & 0 & 1
\end{array}\right).$$
By Theorem \ref{class2}, it remains to multiply these matrices.
	
In the case $\alpha_2=\pm 1$ we need to use Proposition \ref{dop} and (\ref{abc1}), (\ref{basis2}).
\end{proof}

The following is immediate from Theorem \ref{mainq} and Propositions  \ref{form}, \ref{dop}:

\begin{theorem}
\label{mainq1}
Let
\begin{equation}
\label{nm1}
w_1=\sqrt{1-\alpha_2^2+(\beta-\alpha_2)^2},\quad\mu_1=\frac{\sin{w_1t}}{w_1},\quad
\nu_1=\frac{1-\cos{w_1t}}{w_1^2}.
\end{equation}
If $\alpha_2\neq\pm 1$ then the geodesic  $\gamma_1(t)=\gamma_1(\alpha_1,\alpha_2,\alpha_3,\beta;t)$, $t\in\mathbb{R}$, of the sub-Riemannian space $(\SO(3)\times\mathbb{R},\rho_1)$  is equal to $(C,v)(t)$, $t\in\mathbb{R}$, where $v(t)=\alpha_2t$ and the columns  $C_j(t)$, $j=1,2,3$, of the matrix $C(t)\in\SO(3)$
are given by the formulas
$$C_1(t)=\left(\begin{array}{c}
\left(1-\nu_2(\alpha_3^2+(\beta-\alpha_2)^2)\right)\cos{\beta t}-\left(\alpha_1\alpha_3\nu_2-(\beta-\alpha_2)\mu_2\right)\sin{\beta t} \\
\left(\alpha_1\alpha_3\nu_2+(\beta-\alpha_2)\mu_2\right)\cos{\beta t}-\left(1-\nu_2(\alpha_1^2+(\beta-\alpha_2)^2)\right)\sin{\beta t} \\
\left(\alpha_1(\beta-\alpha_2)\nu_2-\alpha_3\mu_2\right)\cos{\beta t}-\left(\alpha_3(\beta-\alpha_2)\nu_2+\alpha_1\mu_2\right)\sin{\beta t} 
\end{array}\right),$$
$$C_2(t)=\left(\begin{array}{c}
\left(1-\nu_2(\alpha_3^2+(\beta-\alpha_2)^2)\right)\sin{\beta t}+\left(\alpha_1\alpha_3\nu_2-(\beta-\alpha_2)\mu_2\right)\cos{\beta t} \\
\left(\alpha_1\alpha_3\nu_2+(\beta-\alpha_2)\mu_2\right)\sin{\beta t}+\left(1-\nu_2(\alpha_1^2+(\beta-\alpha_2)^2)\right)\cos{\beta t} \\
\left(\alpha_1(\beta-\alpha_2)\nu_2-\alpha_3\mu_2\right)\sin{\beta t}+\left(\alpha_3(\beta-\alpha_2)\nu_2+\alpha_1\mu_2\right)\cos{\beta t} 
\end{array}\right),$$
$$C_3(t)=\left(\begin{array}{c}
\alpha_1(\beta-\alpha_2)\nu_2+\alpha_3\mu_2 \\
\alpha_3(\beta-\alpha_2)\nu_2-\alpha_1\mu_2 \\
1-\nu_2(\alpha_1^2+\alpha_3^2) 
\end{array}\right).$$
	
If $\alpha_2=\pm 1$  then the geodesic of the sub-Riemannian space $(\SO(3)\times\mathbb{R},\rho_1)$	is equal to $\gamma_1(t)=(C,v)(t)$, $t\in\mathbb{R}$,  where $v(t)=\alpha_2t$ and
$$C(t)=\left(\begin{array}{ccc}
\cos{\alpha_2 t} & \sin{\alpha_2 t} & 0  \\
-\sin{\alpha_2 t} & \cos{\alpha_2 t} & 0  \\
0 & 0 & 1  
\end{array}\right).$$
\end{theorem}

\section{Cut locus and the first caustic in $\SO(3)\times\mathbb{R}$}

The group $\SU(2)$ is a simply connected covering group of $\SO(3)$; the double covering $\Pi:\SU(2)\rightarrow\SO(3)$ can be defined as follows:
\begin{equation}
\label{Matr}
\Pi(A,B)=\left(\begin{array}{ccc}
A_1^2-A^2_2+B_1^2-B^2_2 & 2(B_1B_2-A_1A_2) & 2(A_1B_2+A_2B_1) \\
2(A_1A_2+B_1B_2) & A_1^2-A^2_2-B^2_1+B_2^2 & 2(A_2B_2-A_1B_1) \\
2(A_2B_1-A_1B_2) & 2(A_1B_1+A_2B_2) & A_1^2+A_2^2-B_1^2-B_2^2
\end{array}\right),
\end{equation}
where $A_1,A_2,B_1,B_2$ are defined by (\ref{AB}).

Then the mapping
\begin{equation}
\label{tildepi}
\tilde{\Pi}:\SU(2)\times\mathbb{R}\rightarrow\SO(3)\times\mathbb{R},\quad \tilde{\Pi}(A,B,v)=(\Pi(A,B),v),
\end{equation}
is a double covering of  $\SO(3)\times\mathbb{R}$ by $\SU(2)\times\mathbb{R}$.

It is easy to check that the differential $d\tilde{\Pi}(\Id)$ is the isomorphism of the Lie algebras $\mathfrak{su}(2)\oplus\mathfrak{\mathbb{R}}$ and
$\mathfrak{so}(3)\oplus\mathfrak{\mathbb{R}}$, translating the basis $E_1,E_2,E_3,E_3$ of  $\mathfrak{su}(2)\oplus\mathfrak{\mathbb{R}}$ given by  (\ref{Basis1}), to the  same name basis of  $\mathfrak{so}(3)\oplus\mathfrak{\mathbb{R}}$ given by  (\ref{basis2}).
Therefore, according to Theorems \ref{class1} and \ref{class2}, the mapping $\tilde{\Pi}:(\SU(2)\times\mathbb{R},d_i)\rightarrow (\SO(3)\times\mathbb{R},\rho_i)$, $i=1,2$, is a local isometry, and
for all $\alpha_1,\alpha_2,\alpha_3,\beta,t\in\mathbb{R}$ such that (\ref{norm}) holds,  we have
\begin{equation}
\label{soot}
\tilde{\Pi}(\tilde{\gamma}_i(\alpha_1,\alpha_2,\alpha_3,\beta;t))=\gamma_i(\alpha_1,\alpha_2,\alpha_3,\beta;t),\quad i=1,2.
\end{equation}

Equality (\ref{soot}) and Corollary \ref{sopr} imply 

\begin{proposition}
	\label{ConjSO3}
	1. The $n$-th conjugate time $t^n_{{\rm conj}}$ along the geodesic $\gamma_i(t)$, $\alpha_2\neq \pm 1$, 
	in $(\SO(3)\times\mathbb{R},\rho_i)$, $i=1,2$, has the form
	$$t^{2m-1}_{\conj}=\frac{2\pi m}{w_i},\quad t^{2m}_{\conj}=\frac{2x_m}{w_i},\quad m\in\mathbb{N},$$
	where $\{x_1,x_2,\dots\}$ are the ascending positive roots of the equation $\tan{x}=x$.
	
	2.  The first conjugate locus  $\Conj_i^{1}$  in $(\SO(3)\times\mathbb{R},d_i)$, $i=1,2$, has the form
	$$\Conj_2^{1}=\{\tilde{\Pi}(A,0,v)\left|\,(A,0,v)\in\SU(2)\times\mathbb{R},\,\,|A|=1,\,\,A\neq 1,\,\,v\in\mathbb{R}\right.\},$$
	$$\Conj_1^{1}=\{\tilde{\Pi}(A,0,v)\left|\,(A,0,v)\in\SU(2)\times\mathbb{R},\,\,|A|=1,\,\,v\notin -2{\rm Arg}(A)\right.\},$$
	where ${\rm Arg}(A)$ is the set of all values of the argument of a complex number  $A$.
\end{proposition}

\begin{proposition}
	\label{sq3}
Let $\gamma_i(t)$, $0\leq t\leq T$, be a shortest arc in $(\SO(3)\times\mathbb{R},\rho_i)$, $i=1,2$, such that
$\gamma_i(0)=\Id$. There exists only one shortest arc $\tilde{\gamma}_i(t)$, $0\leq t\leq T$, in $(\SU(2)\times\mathbb{R},d_i)$ such that $\tilde{\gamma}_i(0)=\Id$ and $\tilde{\Pi}(\tilde{\gamma}_i(t))=\gamma_i(t)$ for every $t\in [0,T]$.
\end{proposition}

\begin{proof}
It follows from (\ref{soot}) that required geodesic $\tilde{\gamma}_i(t)$, $0\leq t\leq T$, exists.
It remains to prove that it is a shortest arc. Let us assume the opposite.
Then there exists a curve  $\tilde{\gamma}_0(t)$, $0\leq t\leq T_0<T$, in $(\SU(2)\times\mathbb{R},d_i)$ such that $\tilde{\gamma}_0(0)=\Id$ and $\tilde{\gamma}_0(T_0)=\tilde{\gamma}_i(T)$. But then the curve $\tilde{\Pi}(\tilde{\gamma}_0(t))$, $0\leq t\leq T_0$, connects $\Id$ and $\gamma_i(T)$, i.e., $d_i(\Id,\gamma_i(T))<T$, a contradiction.
\end{proof}

The following is immediate from Propositions \ref{ssq3} and \ref{sq3}.

\begin{proposition}
\label{ssq4}
Let $\gamma_1(\varphi_0,\alpha_2,\beta;t)$, $0\leq t\leq T$, be a noncontinuable shortest arc in
$(SO(3)\times\mathbb{R},\rho_1)$.
Then $\gamma_2(\varphi_0,\alpha_2,\beta-\alpha_2;t)$, $0\leq t\leq T$, is a noncontinuable shortest arc in 
$(SO(3)\times\mathbb{R},\rho_2)$. The converse is also true.
\end{proposition}

\begin{proposition}
\label{CutSO3}
The cut locus $\Cut_2$ of the sub-Riemannian space $(\SO(3)\times\mathbb{R},\rho_2)$ corresponding
to ${\rm Id}$ is the set
$$\Cut_2=\Cut_2^{\loc}\cup\Cut_2^{\glob},$$
where
$$\Cut_2^{\loc}=\left\{\tilde{\Pi}(A,0,v)\left|\,(A,0,v)\in\SU(2)\times\mathbb{R},\,\,|A|=1,\,\,A\neq\pm 1,\,\,v\in\mathbb{R}\right.\right\}$$
$$=\left\{(C,v)\in\SO(3)\times\mathbb{R}\left|\,\,
C=\left(
\begin{array}{ccc}
\cos\psi & -\sin\psi & 0 \\
\sin\psi & \cos\psi & 0 \\
0 & 0 & 1
\end{array}
\right),\,\,\psi\in (0,2\pi),\,v\in\mathbb{R}\right.\right\},$$
$$\Cut_2^{\glob}=\left\{\tilde{\Pi}(A,B,v)\left|\,(A,0,v)\in\SU(2)\times\mathbb{R},\,\,\re(A)=0,\,v\in\mathbb{R}\right.\right\}$$
$$=\left\{(C,v)\in\SO(3)\times\mathbb{R} \left|C=C^{T},\,\Tr(C)=-1,\,v\in\mathbb{R}\right.\right\}.$$
\end{proposition}

\begin{proof} 
Let us prove that $\Cut_2^{\loc}\subset\Cut_2$. Let $g\in\Cut_2^{\loc}$ and $\gamma_2(t)$, $0\leq t\leq T$, be a shortest arc in  $(\SO(3)\times\mathbb{R},\rho_2)$ joining $\Id$ and $g$.
By virtue of Proposition \ref{sq3}, there exists  only one shortest arc $\tilde{\gamma}_2(t)$, $0\leq t\leq T$, of sub-Riemannian space
$(\SU(2)\times\mathbb{R},d_2)$ such that $\tilde{\gamma}_2(0)=\Id$ and $\tilde{\Pi}(\tilde{\gamma}_2(t))=\gamma_2(t)$ for every $t\in [0,T]$. It follows from Proposition  \ref{ConjSO3} and the proof of Proposition  \ref{cut} that $\tilde{g}:=\tilde{\gamma}_2(T)$ belongs to the cut locus (for $\Id$) in $(\SU(2)\times\mathbb{R},d_2)$, and there exists another shortest arc $\tilde{\gamma}^{\ast}_2(t)$, $0\leq t\leq T$, joining
$\Id$ and $\tilde{g}$ such that $\tilde{\Pi}(\tilde{g})=g$.
Thus the geodesic $\tilde{\Pi}(\tilde{\gamma}^{\ast}_2(t))\neq\gamma_2(t)$, $0\leq t\leq T$, in
$(\SO(3)\times\mathbb{R},\rho_2)$ is a shortest arc joining $\Id$ and $g$. Then $g\in\Conj_2^{1}$ by Proposition \ref{bz}.

Let us prove that $\Cut_2^{\glob}\subset\Cut_2$.  Let $g\in\Cut_2^{\glob}$ and $\tilde{\gamma}_2(t)$, $0\leq t\leq T$, be a shortest arc in 
$(\SO(3)\times\mathbb{R},\rho_2)$ joining  $\Id$ and $g$.
By virtue of Proposition \ref{sq3}, there exists  only one shortest arc  $\tilde{\gamma}_2(t)$, $0\leq t\leq T$, in $(\SU(2)\times\mathbb{R},d_2)$ such that $\tilde{\gamma}_2(0)=\Id$ and $\tilde{\Pi}(\tilde{\gamma}_2(t))=\gamma_2(t)$ for every $t\in [0,T]$. Put
$\tilde{\gamma}_2(T)=(A,B,v)$. It follows from the definition of $\Cut_2^{\glob}$ that $\re(A)=0$. If $B=0$ then $\tilde{\Pi}(\tilde{\gamma}_2(T))\in\Cut_2^{\loc}$ or $g=\Id$; so we can assume that $B\neq 0$. By (\ref{form3}) and Remark \ref{zam0}, we can assume that
$\tilde{\gamma}_2(t)=\tilde{\gamma}_i(\varphi_0,\alpha_2,\beta;t)$.
From (\ref{fform1}) and (\ref{fform2}) we obtain that
$$\tilde{\gamma}_2(\pi+\varphi_0+\beta T,\alpha_2,-\beta;T)=(\bar{A},-B,v)=(-A,-B,v).$$
Then by virtue of (\ref{Matr}) and (\ref{tildepi}),
$\tilde{\Pi}(\tilde{\gamma}_2(\pi+\varphi_0+\beta T,\alpha_2,-\beta;T))=\gamma_2(T)=g,$
i.e. the points $\Id$ and $g$ in  $(\SO(3)\times\mathbb{R},\rho_2)$ are connected by two different shortest paths parametrized by
arclength. By Proposition \ref{bz}, $g\in\Cut_2$.

Now let $g\in\SO(3)\times\mathbb{R}$, but $g\notin\Cut_2^{\loc}\cup\Cut_2^{\glob}$.
We assume that in $(\SO(3)\times\mathbb{R},\rho_2)$ there exist two distinct shortest arcs
$\tilde{\Pi}(\tilde{\gamma}_2(\varphi_0,\alpha_2,\beta;t))$ and $\tilde{\Pi}(\tilde{\gamma}^{\ast}_2(\varphi^{\ast}_0,\alpha_2^{\ast},\beta^{\ast};t))$, $0\leq t\leq T$, connecting $\Id$ and $g$.
Let $\tilde{\gamma}_2(\varphi_0,\alpha_2,\beta;T):=(A,B,v)$ and $\tilde{\gamma}^{\ast}_2(\varphi^{\ast}_0,\alpha_2^{\ast},\beta^{\ast};T):=(A^{\ast},B^{\ast},v^{\ast})$.
By (\ref{Matr}), we have $(A,B,v)=(A^{\ast},B^{\ast},v^{\ast})$ or $(A,B,v)=(-A^{\ast},-B^{\ast},v^{\ast})$. In the first case
$\tilde{\gamma}_2(\varphi_0,\alpha_2,\beta;T)$ belongs to the cut locus of the sub-Riemannian space  $(\SU(2)\times\mathbb{R},d_2)$, thus $g\in\Cut_2^{\loc}$, which is impossible. In the second case without loss of generality, we can assume that $\re(A)<0$. Since $\tilde{\gamma}_2(\varphi_0,\alpha_2,\beta;t)$, $0\leq t\leq T$,
connects $\Id$ and $(A,B,v)$, and $\re(A)<0$ then there exists $t_1\in (0,T)$ such that
$\re (A_1)=0$, where $\tilde{\gamma}_2(\varphi_0,\alpha_2,\beta;t_1)=(A_1,B_1,v_1)$. Then $\tilde{\Pi}(\tilde{\gamma}_2(\varphi_0,\alpha_2,\beta;t_1))\in\Cut_2^{\glob}$. This contradicts the fact that $\tilde{\Pi}(\tilde{\gamma}_2(\varphi_0,\alpha_2,\beta;t))$, $0\leq t\leq T$, is a shortest arc. 
\end{proof}

\begin{proposition}
	\label{CutSO32}
	The cut locus $\Cut_1$ of the sub-Riemannian space $(\SO(3)\times\mathbb{R},\rho_1)$ corresponding
	to ${\rm Id}$ is the set
	$$\Cut_1=\Cut_1^{\loc}\cup\Cut_1^{\glob},$$
	where
	$$\Cut_1^{\loc}=\left\{(C,v)\left|\,\,
	C=\left(
	\begin{array}{ccc}
	\cos\psi & -\sin\psi & 0 \\
	\sin\psi & \cos\psi & 0 \\
	0 & 0 & 1
	\end{array}
	\right),\,\,v\in\mathbb{R},\,\,\psi+v\neq 2\pi n,\,\,n\in\mathbb{N}\right.\right\},$$
$$\Cut_1^{{\rm glob}}=
\left\{(C^{\ast},v)\mid\ C^{\ast}=\left(
\begin{array}{ccc}
	c_{11}\cos{v}-c_{12}\sin{v} & c_{11}\sin{v}+c_{12}\cos{v} & c_{13} \\
	c_{21}\cos{v}-c_{22}\sin{v} & c_{21}\sin{v}+c_{22}\cos{v} & c_{23} \\
	c_{31}\cos{v}-c_{32}\sin{v} & c_{31}\sin{v}+c_{32}\cos{v} & c_{33}
\end{array}\right),\right.$$
$$\left.\quad C=(c_{ij})\in SO(3),\,\,C=C^{T},\,\,{\rm Tr}(C)=-1,\,\,v\in\mathbb{R}
\right\}.$$
\end{proposition}

\begin{proof}
	It follows from the definition of the cut set, Propositions \ref{form}, \ref{ssq4} and Theorem \ref{mainn} that $g_2=(C,v)\in {\rm Cut}_2$ if and only if $g_1=(C,v)\exp(-vE_3)=(C^{\ast},v)\in {\rm Cut}_1$, где
	$$C^{\ast}=\left(
	\begin{array}{ccc}
	c_{11}\cos{v}-c_{12}\sin{v} & c_{11}\sin{v}+c_{12}\cos{v} & c_{13} \\
	c_{21}\cos{v}-c_{22}\sin{v} & c_{21}\sin{v}+c_{22}\cos{v} & c_{23} \\
	c_{31}\cos{v}-c_{32}\sin{v} & c_{31}\sin{v}+c_{32}\cos{v} & c_{33}
	\end{array}
	\right).$$
	From here and Proposition \ref{CutSO3} it follows that $g_2\in{\rm Cut}_2^{{\rm loc}}$ ($g_2\in{\rm Cut}_2^{{\rm glob}}$)  if and only if  $g_1\in {\rm Cut}_1^{{\rm loc}}$ (respectively $g_1\in {\rm Cut}_1^{{\rm glob}}$).

\end{proof}

\section{Sub-Riemannian distance and noncontinuable shortest arcs on the Lie group $\SO(3)\times\mathbb{R}$}

\begin{proposition}
\label{k2}
For every  $(C,v)\in\SO(3)\times\mathbb{R}$,
$$d_1(\Id,(C,v))=d_2(\Id,(\tilde{C},v)),$$
where
$$\tilde{C}=
\left(
\begin{array}{ccc}
c_{11}\cos{v}+c_{12}\sin{v} & c_{12}\cos{v}-c_{11}\sin{v} & c_{13} \\
c_{21}\cos{v}+c_{22}\sin{v} & c_{22}\cos{v}-c_{21}\sin{v} & c_{23} \\
c_{31}\cos{v}+c_{32}\sin{v} & c_{32}\cos{v}-c_{31}\sin{v} & c_{33}
\end{array}
\right).$$
\end{proposition}

\begin{proof}
By analogy with the proof of Proposition \ref{k}, we see that	
$d_1(\Id,(C,v))=d_2(\Id,(C,v)\exp(vE_3)),$
where $E_3$ is given by (\ref{basis2}). It remains to note that \linebreak
$(C,v)\exp(vE_3)=(\tilde{C},v)$.
\end{proof}

There is a consequence of Theorem \ref{os}:

\begin{theorem}
\label{ppq}
$d_2^2(\Id,(C,v))=v^2+d_2^2(\Id,(C,0))$ for all $(C,v)\in\SO(3)\times\mathbb{R}$.
\end{theorem}

The following is immediate from here and \cite[Theorem~4.29]{Sachkov}:

\begin{theorem}
\label{dist}
Let $g=(C,v)\in SO(3)\times\mathbb{R}$,  $g\neq\Id$.
Then
	
1. If $c_{33}=-1$ then $d_2(\Id,g)=\sqrt{v^2+\pi^2}$.
	
2. If $c_{33}=1$ then
$d_2(\Id,g)=\sqrt{v^2+\frac{4\pi^2}{1+\xi^2}},$
where $\xi$ is the unique solution for system of equations
$$\cos\frac{\pi\xi}{\sqrt{1+\xi^2}}=-\frac{1}{2}\sqrt{1+c_{11}+c_{22}+c_{33}},$$
$$\sin\frac{\pi\xi}{\sqrt{1+\xi^2}}=\frac{1}{2}\sgn(c_{21}-c_{12})\sqrt{1-c_{11}-c_{22}+c_{33}}.$$
	
3. If $-1<c_{33}<1$ and $\cos\left(\pi\sqrt{(1+c_{33})/2}\right)=-(c_{11}+c_{22})/(1+c_{33})$ then 
$$d_2(\Id,g)=\sqrt{v^2+\pi^2(1-c_{33})/2}.$$
	
4. If $-1<c_{33}<1$ and $\cos\left(\pi\sqrt{(1+c_{33})/2}\right)>-(c_{11}+c_{22})/(1+c_{33})$ then 
$$d_2(\Id,g)=\sqrt{v^2+\frac{4}{1+\xi^2}\left(\sin^{-1}\sqrt{\frac{1}{2}(1-c_{33})(1+\xi^2)}\right)^2},$$
where $\xi$ is unique solution for system of equations
$$\cos{\left(-\frac{\xi}{\sqrt{1+\xi^2}}\sin^{-1}\sqrt{\frac{1}{2}(1-c_{33})(1+\xi^2)}+
\sin^{-1}\left(\xi\sqrt{\frac{1-c_{33}}{1+c_{33}}}\right)\right)}$$
$$=\sqrt{\frac{1+c_{11}+c_{22}+c_{33}}{2(1+c_{33})}},$$
$$\sin{\left(-\frac{\xi}{\sqrt{1+\xi^2}}\sin^{-1}\sqrt{\frac{1}{2}(1-c_{33})(1+\xi^2)}+
\sin^{-1}\left(\xi\sqrt{\frac{1-c_{33}}{1+c_{33}}}\right)\right)}$$
$$\sgn(c_{21}-c_{12})\sqrt{\frac{1-c_{11}-c_{22}+c_{33}}{2(1+c_{33})}}.$$
	
5. If $-1<c_{33}<1$ and $\cos\left(\pi\sqrt{(1+c_{33})/2}\right)<-(c_{11}+c_{22})/(1+c_{33})$ then 
$$d_2(\Id,g)=\sqrt{v^2+\frac{4}{1+\xi^2}\left(\pi-\sin^{-1}\sqrt{\frac{1}{2}(1-c_{33})(1+\xi^2)}\right)^2},$$
where $\xi$ is the unique solution for system of equations
$$\cos{\left(\frac{\xi}{\sqrt{1+\xi^2}}\left(\pi-\sin^{-1}\sqrt{\frac{1}{2}(1-c_{33})(1+\xi^2)}\right)+
\sin^{-1}\left(\xi\sqrt{\frac{1-c_{33}}{1+c_{33}}}\right)\right)}$$
$$=-\sqrt{\frac{1+c_{11}+c_{22}+c_{33}}{2(1+c_{33})}},$$
$$\sin{\left(\frac{\xi}{\sqrt{1+\xi^2}}\left(\pi-\sin^{-1}\sqrt{\frac{1}{2}(1-c_{33})(1+\xi^2)}\right)+
\sin^{-1}\left(\xi\sqrt{\frac{1-c_{33}}{1+c_{33}}}\right)\right)}$$
$$=\sgn(c_{21}-c_{12})\sqrt{\frac{1-c_{11}-c_{22}+c_{33}}{2(1+c_{33})}}.$$
\end{theorem}

Using Proposition \ref{k2} and Theorem \ref{dist}, we can obtain some exact formulas for $d_1(\Id,(C,v))$.
We will not present them because of their bulkiness.

By Proposition \ref{dop2}, to find shortest arcs in $(\SO(3)\times\mathbb{R},d_i)$, $i=1,2$, 
it suffices to study the  geodesic segments
$\gamma_i(t)$, $0\leq t\leq T$, where $\gamma_i(0)=\Id$.

The following proposition is proved similarly to item~1 of Proposition \ref{cut}.

\begin{proposition}
\label{ss}
If $\alpha_2=\pm 1$ then each geodesic $\gamma_i(t)=\exp(t\alpha_2e_2)$, $t\in\mathbb{R}$, is a metric line of the  
sub-Riemannian space $(\SO(3)\times\mathbb{R},\rho_i)$, $i=1,2$, i.e., for all $t_0,t_1\in\mathbb{R}$ the geodesic segment $\gamma_i(t)=\exp(t\alpha_2e_2)$, $t_0\leq t\leq t_1$,  is a  shortest arc.	
\end{proposition}

\begin{theorem}
\label{kr}
Let $\alpha_2\neq\pm 1$ and $\gamma_2(t)=\gamma_2(\varphi_0,\alpha_2,\beta;t)$, $0\leq t\leq T$, be a noncontinuable shortest arc of the sub-Riemannian space $(\SO(3)\times\mathbb{R},\rho_2)$. Then

1. If $|\beta|\geq\sqrt{\frac{1-\alpha_2^2}{3}}$ then $T=\frac{2\pi}{\sqrt{1-\alpha_2^2+\beta^2}}$.

2. If $|\beta|<\sqrt{\frac{1-\alpha_2^2}{3}}$ then $T\in \left(0,\frac{2\pi}{\sqrt{1-\alpha_2^2+\beta^2}}\right)$ is unique solution of the equation
\begin{equation}
\label{T2}
\cos\frac{T\sqrt{1-\alpha_2^2+\beta^2}}{2}\cos\frac{|\beta|T}{2}+\frac{|\beta|}{\sqrt{1-\alpha_2^2+\beta^2}}\sin\frac{T\sqrt{1-\alpha_2^2+\beta^2}}{2}\sin\frac{|\beta|T}{2}=0.
\end{equation}
3. $T=T(|\beta|)$ is continuous function increasing if $0\leq |\beta|\leq\sqrt{\frac{1-\alpha_2^2}{3}}$ and decreasing if 
$|\beta|\geq\sqrt{\frac{1-\alpha_2^2}{3}}$.
\end{theorem}

\begin{proof}
Let $\alpha_2\neq\pm 1$ and $\gamma_2(t)=\gamma_2(\alpha_1,\alpha_2,\alpha_3,\beta;t)$, $0\leq t\leq T$, be a noncontinuable shortest arc. By Proposition \ref{CutSO3}, $g:=\gamma_2(T)$ belongs to the union of $\Cut_2^{\loc}$ and $\Cut_2^{\glob}$.
Then, by Theorem \ref{mainn}, we get: 

1) if $g\in\Cut_2^{\loc}$ then $m_2(T)=0$ and $T=\frac{2\pi}{\sqrt{1-\alpha_2^2+\beta^2}}$; 

2) if $g\in\Cut_2^{\glob}$ then  (\ref{T2}) holds. 

Consequently, $T=\min\{T_1,T_2\}$, where $T_1=\frac{2\pi}{\sqrt{1-\alpha_2^2+\beta^2}}$ and $T_2$ is the smallest positive root of (\ref{T2}).

First of all, note that for $\beta=0$ equality (\ref{T2}) can be written as
$\cos\frac{T\sqrt{1-\alpha_2^2}}{2}=0$; i.e., $T_2=\frac{\pi}{\sqrt{1-\alpha_2^2}}<T_1$ and $T=T_2$.

We need the following lemma:

\begin{lemma}
\label{aa}
Let $\alpha_2\neq\pm 1$ and $\beta\neq 0$.
Define the function $F(t)$ on the segment  $\left[0,\frac{2\pi}{\sqrt{1-\alpha_2^2+\beta^2}}\right]$ by the formula
$$F(t)=\cos\frac{t\sqrt{1-\alpha_2^2+\beta^2}}{2}\cos\frac{|\beta|t}{2}+\frac{|\beta|}{\sqrt{1-\alpha_2^2+\beta^2}}\sin\frac{t\sqrt{1-\alpha_2^2+\beta^2}}{2}\sin\frac{|\beta|t}{2}.$$
Then

1. If $|\beta|>\sqrt{\frac{1-\alpha_2^2}{3}}$ then $F(t)>0$ on its whole domain.

2. If $|\beta|=\sqrt{\frac{1-\alpha_2^2}{3}}$ then $t=\frac{2\pi}{\sqrt{1-\alpha_2^2+\beta^2}}$ is the unique root of $F(t)$.

3. If $|\beta|<\sqrt{\frac{1-\alpha_2^2}{3}}$ then $F(t)$ has the unique root on the interval $\left(0,\frac{2\pi}{\sqrt{1-\alpha_2^2+\beta^2}}\right)$.
\end{lemma}

\begin{proof}
It is easy to see that
$$F^{\prime}(t)=-\frac{(1-\alpha_2^2)}{2\sqrt{1-\alpha_2^2+\beta^2}}\sin\frac{t\sqrt{1-\alpha_2^2+\beta^2}}{2}\cos\frac{|\beta|t}{2}.$$

If $|\beta|>\sqrt{\frac{1-\alpha_2^2}{3}}$ then $\frac{\pi}{|\beta|}<\frac{2\pi}{\sqrt{1-\alpha_2^2+\beta^2}}$, and $F(t)$
decreases on the segment $\left[0,\frac{\pi}{|\beta|}\right]$ and increases on the segment $\left[\frac{\pi}{|\beta|},\frac{2\pi}{\sqrt{1-\alpha_2^2+\beta^2}}\right]$. Moreover, $F(0)=1$ and 
$$F\left(\frac{\pi}{|\beta|}\right)=\frac{|\beta|}{\sqrt{1-\alpha_2^2+\beta^2}}\sin\frac{\pi\sqrt{1-\alpha_2^2+\beta^2}}{2|\beta|}>0.$$ 
This implies item~1 of  Lemma \ref{aa}.

If $|\beta|=\sqrt{\frac{1-\alpha_2^2}{3}}$ then $\frac{\pi}{|\beta|}=\frac{2\pi}{\sqrt{1-\alpha_2^2+\beta^2}}$,   and the function $F(t)$
decreases on the segment $\left[0,\frac{2\pi}{\sqrt{1-\alpha_2^2+\beta^2}}\right]$. Moreover,  $F\left(\frac{2\pi}{\sqrt{1-\alpha_2^2+\beta^2}}\right)=0$.
This implies item~2 of  Lemma \ref{aa}.
 
If $|\beta|<\sqrt{\frac{1-\alpha_2^2}{3}}$ then $\frac{\pi}{|\beta|}>\frac{2\pi}{\sqrt{1-\alpha_2^2+\beta^2}}$, and the function   $F(t)$
decreases on the segment $\left[0,\frac{2\pi}{\sqrt{1-\alpha_2^2+\beta^2}}\right]$. Moreover, $F(0)=1$ and 
$$F\left(\frac{2\pi}{\sqrt{1-\alpha_2^2+\beta^2}}\right)=-\cos\frac{\pi|\beta|}{\sqrt{1-\alpha_2^2+\beta^2}}<0.$$ 
This implies item~3 of  Lemma \ref{aa}.
\end{proof}

It follows from (\ref{T2}) and Lemma \ref{aa} that items~1 and 2 of Theorem  \ref{kr} and the continuity of the function $T(|\beta|)$ hold on the ray $[0,+\infty)$.

Because of item~1 of Theorem  \ref{kr}, the $T(|\beta|)$ decreases if $|\beta|\geq\sqrt{\frac{1-\alpha_2^2}{3}}$. 
To differentiate 
(\ref{T2}) with respect to $|\beta|$ and  express $T^{\prime}(|\beta|)$, using (\ref{T2}) again, we get
$$T^{\prime}(|\beta|)=\frac{{\rm tg}(\beta |T|/2)}{1-\alpha_2^2+\beta^2}\left(2-T\sqrt{1-\alpha_2^2+\beta^2}{\rm ctg}\left(\frac{T}{2}\sqrt{1-\alpha_2^2+\beta^2}\right)\right)=$$
$$\frac{2{\rm tg}(|\beta| T/2)}{1-\alpha_2^2+\beta^2}\left(1+\frac{\beta|T|}{2}{\rm tg}\frac{\beta |T|}{2}\right)>0$$
for $0<|\beta|<\sqrt{\frac{1-\alpha_2^2}{3}}$, because then  $0<\frac{|\beta|T}{2}<\frac{\pi}{2}$.

Theorem \ref{kr} is proved.
\end{proof}

By virtue of Proposition \ref{ssq4}, if  $\gamma_1(t)=\gamma_2(\varphi_0,\alpha_2,\beta;t)$, $0\leq t\leq T$, is a noncontinuable shortest
arc of the sub-Riemannian space $(\SO(3)\times\mathbb{R},\rho_1)$, then $T$ satisfies items~1--3 of Theorem \ref{kr} after replacing $\beta$ by $\beta-\alpha_2$.

\vspace{2mm}

{\bf Acknowledgments.}
The author thanks Professor V.N.~Berestovskii and 
\linebreak E.A.~Meshcheryakov for useful discussions and remarks.

\end{document}